\documentclass[final]{amsart}
\usepackage{amssymb}
\usepackage{enumerate}

\makeatletter%
\def\th@mytheorem{%
 \let\thm@indent\noindent
 \thm@headfont{\bfseries}
  \itshape
}
\def\th@myremark{%
 \let\thm@indent\noindent
 \thm@headfont{\bfseries}
}
\makeatother%
\theoremstyle{mytheorem}
\newtheorem{Theorem}{Theorem}[section]
\newtheorem{Lemma}[Theorem]{Lemma}

\theoremstyle{myremark}
\newtheorem{Remark}[Theorem]{Remark}
\def\bs{\boldsymbol}
\def\bu{$\bullet$\quad}
\def\BX{\bs{X}}
\def\BY{\bs{Y}}
\def\BZ{\bs{Z}}

\def\CO{\mathcal O}
\def\CC{\mathbb C}
\def\PP{\mathbb P}
\DeclareSymbolFont{EulerScriptBold}{U}{eus}{b}{n}
\DeclareSymbolFontAlphabet\eusb{EulerScriptBold}
\def\EX{\eusb X}

\def\halfskip{\vskip 6pt plus 1pt minus 1pt}
\def\id{\operatorname{id}}
\def\NN{\mathbb N}
\def\XX{\mathbb X}
\def\Omega{\varOmega}
\let\phi=\varphi
\def\PSH{\mathcal{PSH}}
\def\ttimes{\times\dots\times}
\def\too{\longrightarrow}
\def\tuu{\longmapsto}
\def\wdht{\widehat}
\def\wdtl{\widetilde}

\def\eps{\varepsilon}
\def\Delta{\varDelta}
\def\Phi{\varPhi}

\begin{document}
\title
{A new cross theorem for separately holomorphic functions}

\author{Marek Jarnicki}
\address{Jagiellonian University, Institute of Mathematics,
{\L}ojasiewicza 6, 30-348 Krak\'ow, Poland} \email{Marek.Jarnicki@im.uj.edu.pl}

\author{Peter Pflug}
\address{Carl von Ossietzky Universit\"at Oldenburg, Institut f\"ur Mathematik,
Postfach 2503, D-26111 Oldenburg, Germany}
\email{pflug@mathematik.uni-oldenburg.de}
\thanks{The research was partially supported by
the grant no.~N N201 361436 of the Polish Ministry of Science and Higher Education
and the DFG-grant 436POL113/103/0-2.}

\subjclass{32D15}

\keywords{separately holomorphic function, cross theorem}

\begin{abstract}
We prove a new cross theorem for separately holomorphic functions.
\end{abstract}

\maketitle


\section{Introduction. Main result}
Throughout the paper we will work in the following geometric context ---
details may be found in \cite{JarPfl2007}.

We fix an integer $N\geq2$ and let $D_j$ be a (connected) \textit{Riemann domain}
over $\CC^{n_j}$, $j=1,\dots,N$. Let $\varnothing\neq A_j\subset D_j$ be
\textit{locally pluriregular}, $j=1,\dots,N$.

We will use the following conventions: $A'_j:=A_1\ttimes A_{j-1}$, $j=2,\dots,N$,
$A''_j:=A_{j+1}\ttimes A_N$, $j=1,\dots,N-1$. Analogously, a point
$a=(a_1,\dots,a_N)\in D_1\ttimes D_N$ may be written as $a=(a'_j,a_j,a''_j)$, where
$a'_j:=(a_1,\dots,a_{j-1})$, $a''_j:=(a_{j+1},\dots,a_N)$ (with obvious exceptions for
$j\in\{1,N\}$).

We define an \textit{$N$--fold cross}
$$
\BX=\BX((D_j,A_j)_{j=1}^N):=\bigcup_{j=1}^N A'_j\times D_j\times A''_j.
$$
One may easily prove that $\BX$ is connected.

We say that a function $f:\BX\too\CC$ is \textit{separately holomorphic on $\BX$}
(we write $f\in\CO_s(\BX)$) if for any $j\in\{1,\dots,N\}$ and $(a'_j,a''_j)\in A'_j\times A''_j$,
the function $D_j\ni z_j\tuu f(a'_j,z_j,a''_j)\in\CC$ is holomorphic in $D_j$.

\halfskip

Let $h_{A_j,D_j}$ denote the relative extremal function of $A_j$
in $D_j$, $j=1,\dots,N$. Recall that $h_{A,D}:=\sup\{u\in\PSH(D):
u\leq1,\;u|_A\leq0\}$, $A\subset D$ (cf. \cite{Kli1991}, \S\;4.5).  Put
$$ \wdht{\BX}:=\{(z_1,\dots,z_N)\in D_1\ttimes D_N:
h^\ast_{A_1,D_1}(z_1)+\dots+h^\ast_{A_N,D_N}(z_N)<1\}, $$ where
${}^\ast$ stands for the upper semicontinuous regularization. One
may prove that $\wdht{\BX}$ is connected and
$\BX\subset\wdht{\BX}$.

\halfskip

The \textit{classical cross theorem} is the following result:

\begin{Theorem}[\cite{Sic1969a}, \cite{Sic1969b}, \cite{Zah1976},
\cite{Sic1981a}, \cite{NguSic1991}, \cite{NguZer1991}, \cite{NguZer1995},
\cite{NTV1997}, \cite{AleZer2001}, \cite{Zer2002}]\label{ThmClassical}
For each $f\in\CO_s(\BX)$ there exists exactly one $\wdht f\in\CO(\wdht{\BX})$ such
that $\wdht f=f$ on $\BX$ and $\sup_{\wdht{\BX}}|\wdht f|=\sup_{\BX}|f|$.
\end{Theorem}

The aim of this note is to extend the above theorem to a class of more general objects, namely
\textit{$(N,k)$--crosses $\BX_{N,k}$} defined for $k\in\{1,\dots,N\}$ as follows:
$$
\BX_{N,k}=\XX_{N,k}((A_j,D_j)_{j=1}^N):=
\bigcup_{\substack{\alpha_1,\dots,\alpha_N\in\{0,1\}\\
\alpha_1+\dots+\alpha_N=k}}\EX_\alpha,
$$
where
$$
\EX_\alpha:=\EX_{1,\alpha_1}\ttimes\EX_{N,\alpha_N},\quad\EX_{j,\alpha_j}:=
\begin{cases} D_j, & \text{ if } \alpha_j=1 \\ A_j, & \text{ if } \alpha_j=0\end{cases}.
$$

Notice that $N$--fold crosses are just $(N,1)$--crosses in the above terminology. Obviously,
$\BX_{N,N}=D_1\ttimes D_N$. Thus, if $N=2$, then in fact we have
only $\BX_{2,1}$.

Recall that the theory of extension of separately holomorphic functions had been first developed
for $N=2$. Then the $N$--fold case (obtained via induction) was considered as a natural
generalization of $\BX_{2,1}$. In our opinion, each of the crosses $\BX_{N,k}$ may be considered
as a \textit{natural} generalization of $\BX_{2,1}$. Consequently, one should try to find an
analogous of the cross theorem for all $(N,k)$--crosses.

We say that a function $f:\BX_{N,k}\too\CC$ is \textit{separately
holomorphic} ($f\in\CO_s(\BX_{N,k})$) if for all $a=(a_1,\dots,a_N)\in A_1\ttimes A_N$ and
$\alpha=(\alpha_1,\dots,\alpha_N)\in\{0,1\}^N$ with $|\alpha|=k$,
the function
$$
D^\alpha:=\prod_{\substack{j\in\{1,\dots,N\}:\\
\alpha_j=1}}D_j\ni z\tuu f(i_{a,\alpha}(z)) $$ is holomorphic,
where $i_{a,\alpha}:D^\alpha\too\EX_\alpha$, $$
i_{a,\alpha}(z):=(w_1,\dots,w_N),\quad w_j:=\begin{cases} z_j, &
\text{ if } \alpha_j=1 \\ a_j, & \text{ if }
\alpha_j=0\end{cases}.
$$

Put
\begin{multline*}
\wdht{\BX}_{N,k}=\wdht{\XX}_{N,k}((A_j,D_j)_{j=1}^N):\\
=\Big\{(z_1,\dots,z_N)\in D_1\ttimes D_N:
\sum_{j=1}^N h_{A_j,D_j}^\ast(z_j)<k\Big\}.
\end{multline*}
Note that $\wdht{\BX}_{N,N}=D_1\ttimes D_N$.

\halfskip

Let $\phi_j:D_j\too\wdtl D_j$ be the envelope of holomorphy
(cf.~\cite{JarPfl2000}, Definition 1.8.1).
Observe that since $\phi_j$ is locally biholomorphic, the set
$\wdtl A_j:=\phi_j(A_j)\subset\wdtl D_j$ is locally pluriregular, $j=1,\dots,N$. Let
$$
\wdtl{\BX}_{N,k}:=\XX_{N,k}((\wdtl A_j,\wdtl D_j)_{j=1}^N),\quad
\wdht{\wdtl\BX}_{N,k}:=\wdht{\XX}_{N,k}((\wdtl A_j,\wdtl D_j)_{j=1}^N).
$$
Put
$$
\phi:D_1\ttimes D_N\too\wdtl D_1\ttimes\wdtl D_N,\quad \phi(z_1,\dots,z_N):=
(\phi_1(z_1),\dots,\phi_N(z_N)).
$$
Note that:

\bu $\phi(\BX_{N,k})\subset\wdtl\BX_{N,k}$,

\bu $\phi(\wdht{\BX}_{N,k})\subset\wdht{\wdtl{\BX}}_{N,k}$ (because
$h_{\wdtl A_j,\wdtl D_j}^\ast\circ\phi_j\leq h_{A_j,D_j}^\ast$, $j=1,\dots,N$).

\halfskip

Our main result is the following \textit{cross theorem for $(N,k)$--crosses}.

\begin{Theorem}\label{ThmNew}
For every $f\in\CO_s(\BX_{N,k})$ there exists exactly one $\wdht f\in\CO(\wdht{\wdtl\BX}_{N,k})$ such that
$\wdht f\circ\phi=f$ on $\BX_{N,k}$ and $\sup_{\wdht{\wdtl\BX}_{N,k}}|\wdht f|=\sup_{\BX_{N,k}}|f|$.
\end{Theorem}

The proof will be presented in \S\;\ref{SectionProofs} and will be
based on Theorem \ref{ThmClassical} and the following technical lemmas
(which might be also useful in other applications).

\begin{Lemma}\label{LemFormula1}
Let $G$ be a Riemann domain over $\CC^n$, let $D\subset\subset G$ be a Riemann domain of holomorphy,
and let $A\subset D$ be non-pluripolar. Put
$$
\Delta(\mu):=\{z\in D: h_{A,D}^\ast(z)<\mu\},\quad 0<\mu\leq1.
$$
Then
$$
h_{\Delta(r),\Delta(s)}^\ast=\max\Big\{0,\frac{h_{A,D}^\ast-r}{s-r}\Big\} \text{ on }
\Delta(s),\quad 0<r<s\leq1.
$$
\end{Lemma}

\begin{Lemma}\label{LemFormula}
Assume additionally that $D_1,\dots,D_N$ are Riemann domains of holomorphy. Then
\begin{multline*}
h_{\wdht{\BX}_{N,k-1}, \wdht{\BX}_{N,k}}^\ast(z)=
\max\Big\{0,\sum_{j=1}^Nh_{A_j,D_j}^\ast(z_j)-k+1\Big\},\\
z=(z_1,\dots,z_N)\in\wdht{\BX}_{N,k},\quad k\in\{2,\dots,N\}.
\end{multline*}
\end{Lemma}

We do not know whether Lemmas \ref{LemFormula1}, \ref{LemFormula} are true for arbitrary
Riemann domains.

\section{Basic properties of $(N,k)$--crosses}
\begin{Remark}\label{RemOtCr}
\begin{enumerate}[(a)]
\item $A_1\ttimes A_N\subset\BX_{N,k}\subset\wdht{\BX}_{N,k}$.

\item $\BX_{N,k-1}\subset\BX_{N,k}$, $\wdht{\BX}_{N,k-1}\subset\wdht{\BX}_{N,k}$, $k=2,\dots,N$.

\item $\BX_{N,k}=(\BX_{N-1,k-1}\times D_N)\cup(\BX_{N-1,k}\times A_N)$, $k=2,\dots,N-1$, $N\geq3$.

\item $\BX_{N,k}$  and $\wdht{\BX}_{N,k}$ are connected.

\item \label{RemOtCrAppr}
If $(D_{j,k})_{k=1}^\infty$ is a sequence of subdomains of $D_j$ such that
$D_{j,k}\nearrow D_j$, $D_{j,k}\supset A_{j,k}\nearrow A_j$, $j=1,\dots,N$,
then $\XX_{N,k}((A_{j,k},D_{j,k})_{j=1}^N)\nearrow\BX_{N,k}$ and
$$
\wdht{\XX}_{N,k}((A_{j,k},D_{j,k})_{j=1}^N)\nearrow\wdht{\BX}_{N,k}.
$$

\item If $D_1,\dots,D_N$ are domains of holomorphy, then $\wdht{\BX}_{N,k}$ is a domain
of holomorphy.
\end{enumerate}
\end{Remark}

\section{Proof of Lemma \ref{LemFormula1}}
Let
$$
L:=h_{\Delta(r),\Delta(s)}^\ast,\quad R:=\max\Big\{0,\frac{h_{A,D}^\ast-r}{s-r}\Big\}.
$$
Put $\Delta[r]:=\{z\in D: h_{A,D}^\ast(z)\leq r\}$. It is clear that
\begin{gather*}
L=h_{\Delta(r),\Delta(s)}\geq h_{\Delta[r],\Delta(s)}^\ast\geq h_{\Delta[r],\Delta(s)}\geq R,\tag{*}\\
L=R=0 \text{ on }\Delta(r),\quad R=0 \text{ on }\Delta[r].
\end{gather*}

Step 1. Reduction to the case $s=1$.

Suppose that $0<r<s<1$. Observe that $\Delta(s)$ is a Riemann region of holomorphy.
Moreover, $h_{A\cap\Delta(s),\Delta(s)}^\ast=(1/s)h_{A,D}^\ast$ on $\Delta(s)$.

Indeed, it obvious that $h_{A\cap\Delta(s),\Delta(s)}^\ast\geq(1/s)h_{A,D}^\ast$ on $\Delta(s)$.
Let $u\in\PSH(\Delta(s))$, $u\leq1$, $u\leq0$ on $A\cap\Delta(s)$. Observe that for every
$z_0\in D\cap\partial(\Delta(s))$ we have $\limsup_{z\to z_0}u(z)\leq1\leq(1/s)h_{A,D}^\ast(z_0)$.
Thus, the function
$$
v:=\begin{cases}\max\{su,h_{A,D}^\ast\} & \text{ on } \Delta(s)\\
h_{A,D}^\ast & \text{ on } D\setminus\Delta(s)\end{cases}
$$
is plurisubharmonic on $D$. It is known that
there exists a pluripolar set $P\subset A$ such that
$h_{A,D}^\ast=0$ on $A\setminus P$ (cf.~\cite{Kli1991}).
Hence, $A\setminus P\subset\Delta(s)$,
$v\leq h_{A\setminus P,D}^\ast=h_{A,D}^\ast$, and therefore,
$h_{A\cap\Delta(s),\Delta(s)}^\ast\leq(1/s)h_{A,D}^\ast$ on $\Delta(s)$.

In particular, $A\cap S$ is not pluripolar for every connected component $S$ of $\Delta(s)$. Hence,
$$
L=h_{\Delta(r),\Delta(s)}=h_{\{h_{A,\Delta(s)}^\ast<r/s\}, \Delta(s)},\quad
R=\max\Big\{0,\frac{h_{A,\Delta(s)}^\ast-r/s}{1-r/s}\Big\}.
$$
Thus the problem for $(D,A,r,s)$ reduces to $(S,A\cap S,r/s,1)$, where $S$ is a connected
component of $\Delta(s)$.

\textit{From now on we assume that $s=1$.}

\halfskip

Step 2. Approximation. Let $A_\nu\nearrow A$, $D_\nu\nearrow D$, where $A_\nu\subset D_\nu$
is non-pluripolar, $\nu\in\NN$.
Suppose that the formula holds for each $(D_\nu,A_\nu,r)$. Then it holds for $(D,A,r)$.

Indeed, we know  that $h_{A_\nu, D_\nu}^\ast\searrow h_{A,D}^\ast$.
Hence $\{h_{A_\nu,D}^\ast<r\}\nearrow\Delta(r)$. Thus
$h_{\{h_{A_\nu,D}^\ast<r\},D}^\ast\searrow h_{\Delta(r),D}^\ast$.

\halfskip

Step 3. The case where $D$ is hyperconvex, $A$ is compact, and $h_{A,D}^\ast$ is continuous.

Let $u\in\PSH(D)$, $u\leq1$, $u\leq0$ on $\Delta[r]$.
Using continuity of $h_{A,D}^\ast$ and \cite{Kli1991}, Proposition 4.5.2,
we easily conclude that $\Delta[r]$ is compact. Let $U:=D\setminus\Delta[r]$.
Observe that for $z_0\in\partial U$ we get
$$
\liminf_{U\ni z\to z_0}(h_{A,D}^\ast(z)-(1-r)u(z)-r)\geq0.
$$
Hence, by the domination principle (cf.~\cite{Kli1991}, Corollary 3.7.4),
$(1-r)u+r\leq h_{A,D}^\ast$ in $U$.
This shows that $h_{\Delta[r],D}\leq R$. Thus, by (*), we get $h_{\Delta[r],D}^\ast\equiv R$
for all $0<r<1$. Observe that $\Delta[r_\nu]\nearrow\Delta(r)$ for $0<r_\nu\nearrow r$.
Consequently, $L\equiv R$.

\halfskip

Step 4. The case where $D$ is hyperconvex and $A$ is compact.

Let $A^{(\eps)}:=\bigcup_{a\in A}\overline{\wdht\PP(a,\eps)}$, where $\wdht\PP(a,\eps)$ stands for
the ``polydisc'' in the sense of the Riemann domain $D$ ($A^{(\eps)}$ is defined for small
$\eps>0$).
By \cite{Kli1991}, Corollary 4.5.9, we know that $h_{A^{(\eps)},D}
=h_{A^{(\eps)},D}^\ast$ is continuous. Thus, using Step 3 and (*), we have
$$
h_{\{h_{A^{(\eps)},D}\leq r\},D}=\max\Big\{0,\frac{h_{A^{(\eps)},D}-r}{1-r}\Big\},
\quad 0<\eps\ll1.
$$
By \cite{Kli1991}, Proposition 4.5.10,
we have $h_{A^{(\eps)},D}\nearrow h_{A,D}$ as $\eps\searrow0$. In particular,
$$
\{h_{A^{(\eps)},D}\leq r\}\searrow\{h_{A,D}\leq r\} \text{  as } \eps\searrow0.
$$
Hence, once again by \cite{Kli1991}, Proposition 4.5.10,
$$
h_{\{h_{A^{(\eps)},D}\leq r\},D}\nearrow h_{\{h_{A^{(\eps)},D}\leq r\},D} \text{ as } \eps\searrow0.
$$
Consequently,
$$
h_{\{h_{A,D}\leq r\},D}=\max\Big\{0,\frac{h_{A,D}-r}{1-r}\Big\}\leq R.
$$
Thus
$h_{\{h_{A,D}\leq r\},D}^\ast\leq R$.
Observe that the set $\{h_{A,D}\leq r\}\setminus\Delta[r]$ is pluripolar.
Consequently, $h_{\Delta[r],D}^\ast\leq R$. We finish the proof as in Step 3.

\halfskip

Step 5. The case where $A$ is open.

We use Step 4 and approximation (Step 2) with
$A_\nu\nearrow A$, $D_\nu\nearrow D$, where $A_\nu\subset\subset D_\nu$ is compact non-pluripolar
and $D_\nu$ is hyperconvex, $\nu\in\NN$.

\halfskip

Step 6. The case where $D$ is hyperconvex and $A\subset\subset D$ is non-pluripolar.

By Step 5 we get
$$
h_{\{h_{\Delta(\eps),D}^\ast<r\},D}^\ast=
\max\Big\{0,\frac{h_{\Delta(\eps),D}^\ast-r}{1-r}\Big\}, \quad 0<\eps<1.
$$
By \cite{Blo2000}, we get
$$
\frac{h_{A,D}^\ast-\eps}{1-\eps}\leq
h_{\Delta(\eps),D}^\ast\leq h_{A,D}^\ast,
$$
in particular,
$h_{\Delta(\eps),D}^\ast\nearrow h_{A,D}^\ast$ as $\eps\searrow0$. Moreover,
$$
\{h_{\Delta(\eps),D}^\ast<\frac{r-\eps}{1-\eps}\}\subset\Delta(r)\subset\{h_{\Delta(\eps),D}^\ast<r\},
\quad 0<\eps<r.
$$
Consequently,
\begin{multline*}
\max\Big\{0,\frac{h_{\Delta(\eps),D}^\ast-\frac{r-\eps}{1-\eps}}{1-\frac{r-\eps}{1-\eps}}\Big\}
=h_{\{h_{\Delta(\eps),D}^\ast<\frac{r-\eps}{1-\eps}\},D}^\ast\\
\geq h_{\Delta(r),D}^\ast\geq h_{\{h_{\Delta(\eps),D}^\ast<r\},D}^\ast=\max\Big\{0,\frac{h_{\Delta(\eps),D}^\ast-r}{1-r}\Big\},
\quad 0<\eps<r.
\end{multline*}
Letting $\eps\searrow0$, we get the required formula.

\halfskip

Step 7. The general case.

We use Step 6 and approximation (Step 2) with
$A_\nu\nearrow A$, $D_\nu\nearrow D$, where $A_\nu\subset\subset D_\nu$ is
non-pluripolar and $D_\nu$ is hyperconvex, $\nu\in\NN$.

The proof of Lemma \ref{LemFormula1} is completed.

\section{Proof of Lemma \ref{LemFormula}}
By Remark \ref{RemOtCr}(\ref{RemOtCrAppr}), we may assume that
$A_j\subset\subset D_j\subset\subset G_j$, where $G_j$ is a Riemann domain over $\CC^{n_j}$,
$j=1,\dots,N$. Fix $2\leq k\leq N$. Let
$$
h_j:=h_{A_j,D_j}^\ast,\; j=1,\dots,N,\quad
h(z_1,\dots,z_N):=h_1(z_1)+\dots+h_N(z_N).
$$
Let
\begin{multline*}
L_{N,k}:=h_{\wdht{\BX}_{N,k-1}, \wdht{\BX}_{N,k}}^\ast,\quad
R_{N,k}(z)=\max\Big\{0,\sum_{j=1}^Nh_{A_j,D_j}^\ast(z_j)-k+1\Big\},\\
z=(z_1,\dots,z_N)\in\wdht{\BX}_{N,k}.
\end{multline*}
It is clear that $L_{N,k}\geq R_{N,k}$ and $L_{N,k}=R_{N,k}=0$ on $\wdht{\BX}_{N,k-1}$.
Fix an $a=(a_1,\dots,a_N)\in\wdht{\BX}_{N,k}\setminus\wdht{\BX}_{N,k-1}$.
We may assume that $h_1(a_1)\leq\dots\leq h_N(a_N)$. Suppose that $h_1(a_1)=\dots=h_s(a_s)=0$
and $h_{s+1}(a_{s+1}),\dots, h_N(a_N)>0$ for an $s\in\{0,\dots,N\}$. Since $h(a)\geq k-1$, we
see that in fact $s\leq N-k\leq N-2$. In particular, if $N=2$, then $s=0$.

Let $\wdht{\BY}_{N-s,p}=\wdht{\XX}_{N-s,p}((A_j,D_j)_{j=s+1}^N)$, $p\in\{k-1,k\}$.
Observe that
$$
\{a_1,\dots,a_s\}\times\wdht{\BY}_{N-s,p}\subset\wdht{\BX}_{N,p},\quad p\in\{k-1,k\}.
$$
Consequently,
$$
h_{\wdht{\BX}_{N,k-1}, \wdht{\BX}_{N,k}}^\ast(a)\leq
h_{\wdht{\BY}_{N-s,k-1}, \wdht{\BY}_{N-s,k}}^\ast(a_{s+1},\dots,a_N).
$$
Thus, if we know that $L_{N-s,k}(a_{s+1},\dots,a_N)\leq R_{N-s,k}(a_{s+1},\dots,a_N)$,
then
$$
L_{N,k}(a)\leq R_{N-s,k}(a_{s+1},\dots,a_N)=R_{N,k}(a).
$$
This reduces the proof to the case $s=0$, i.e.~$h_j(a_j)>0$, $j=1,\dots,N$.

Put
$$
\Delta_{j,t}:=\{z_j\in D_j: h_j(z_j)<t\},\quad j=1,\dots,N.
$$
Take $0<r_j<s_j\leq1$, $j=1,\dots,N$, such that $r_1+\dots+r_N=k-1$ and
$s_1+\dots+s_N=k$. Observe that
$$
\Delta_{1,r_1}\ttimes\Delta_{N,r_N}\subset\wdht{\BX}_{N,k-1},\quad
\Delta_{1,s_1}\ttimes\Delta_{N,s_N}\subset\wdht{\BX}_{N,k}.
$$
Hence, using the product property for the relative extremal function
(cf.~\cite{Edi2002}, Theorem 4.1) and Lemma \ref{LemFormula1}, we get
\begin{align*}
L_{N,k}(z)&\leq h_{\Delta_{1,r_1}\ttimes\Delta_{N,r_N},\Delta_{1,s_1}\ttimes\Delta_{N,s_N}}^\ast(z)\\
&=\max\{h_{\Delta_{1,r_1},\Delta_{1,r_1}}^\ast(z_1),\dots,h_{\Delta_{N,r_N},\Delta_{N,r_N}}^\ast(z_N)\}\\
&=\max\Big\{0,\frac{h_1(z_1)-r_1}{s_1-r_1},\dots,\frac{h_N(z_N)-r_N}{s_N-r_N}
\Big\},\\
&\hskip150pt z=(z_1,\dots,z_N)\in \Delta_{1,s_1}\ttimes\Delta_{N,s_N}.
\end{align*}
Observe that there exist numbers $s_1,\dots,s_N\in(0,1]$ such that $s_1+\dots+s_N=k$ and
$$
h_j(a_j)<s_j<\frac{h_j(a_j)}{h(a)-k+1},\quad j=1,\dots,N.
$$

Indeed, since the case where $h(a)=k-1$ is trivial, we may assume that $h(a)>k-1$.
Note that $h_j(a_j)<\frac{h_j(a_j)}{h(a)-k+1}$, $j=1,\dots,N$. Suppose that
$$
\frac{h_j(a_j)}{h(a)-k+1}\leq1,\; j=1,\dots,\sigma,\quad
\frac{h_j(a_j)}{h(a)-k+1}>1,\;j=\sigma+1,\dots,N,
$$
for a $\sigma\in\{0,\dots,N\}$.
Observe that
$$
\sum_{j=1}^N\frac{h_j(a_j)}{h(a)-k+1}=\frac{h(a)}{h(a)-k+1}>k,
$$
so the case $\sigma=N$ is simple. Thus, assume that $\sigma\leq N-1$.
We only need do show that
$$
\Big(\sum_{j=1}^\sigma\frac{h_j(a_j)}{h(a)-k+1}\Big)+N-\sigma>k.
$$
The case where $\sigma\leq N-k$ is obvious. Thus assume that $\sigma\geq N-k+1$.
We have to show that
\begin{multline*}
\sum_{j=1}^\sigma h_j(a_j)>(h(a)-k+1)(k-N+\sigma)\\
=(k-1-N+\sigma)h(a)+
\Big(\sum_{j=1}^\sigma h_j(a_j)\Big)+\Big(\sum_{j=\sigma+1}^Nh_j(a_j)\Big)+(-k+1)(k-N+\sigma),
\end{multline*}
or equivalently,
$$
(k-1-N+\sigma)h(a)+\Big(\sum_{j=\sigma+1}^Nh_j(a_j)\Big)<(k-1)(k-N+\sigma).
$$
We have
\begin{multline*}
(k-1-N+\sigma)h(a)+\Big(\sum_{j=\sigma+1}^Nh_j(a_j)\Big)\\
<(k-1-N+\sigma)k+N-\sigma\leq(k-1)(k-N+\sigma),
\end{multline*}
which gives the required inequality.

\halfskip

Now, define
$$
r_j:=\frac{h_j(a_j)-s_j(h(a)-k+1)}{k-h(a)}, \quad j=1,\dots,N.
$$
Then:

\bu $r_j>0$ because $s_j<\frac{h_j(a_j)}{h(a)-k+1}$,

\bu $r_j<s_j$ because $h_j(a)<s_j$,

\bu $r_1+\dots+r_N=k-1$,

\bu $\frac{h_j(a_j)-r_j}{s_j-r_j}=h(a)-k+1$, $j=1,\dots,N$.

Thus
\begin{multline*}
L_{N,k}(a)\leq\max\Big\{0,\frac{h_1(a_1)-r_1}{s_1-r_1},\dots,\frac{h_N(a_N)-r_N}{s_N-r_N}\Big\}\\
=\max\{0,h(a)-k+1\}=R_{N,k}(a).
\end{multline*}
The proof of Lemma \ref{LemFormula} is completed.

\section{Proof of Theorem \ref{ThmNew}}\label{SectionProofs}
First we prove that
for each function $f\in\CO_s(\BX_{N,k})$ there exists exactly one $\wdtl f\in\CO_s(\wdtl\BX_{N,k})$
such that $\wdtl f\circ\phi\equiv f$ and $\sup_{\wdtl\BX_{N,k}}|\wdtl f|=\sup_{\BX_{N,k}}|f|$.

Indeed, fix an $f\in\CO_s(\BX_{N,k})$. Take $a=(a_1,\dots,a_N), b=(b_1,\dots,b_N)\in A_1\ttimes A_N$ and
$\alpha=(\alpha_1,\dots,\alpha_N), \beta=(\beta_1,\dots,\beta_N)\in\{0,1\}^N$ with
$|\alpha|=|\beta|=k$. To simplify notation, suppose that $\alpha=(1,\dots,1,0,\dots,0)$.

Observe that if $\phi_j(a_j)=\phi_j(b_j)$, $j=k+1,\dots,N$, then
$f(\cdot,a_{k+1},\dots,a_N)\equiv f(\cdot,b_{k+1},\dots,b_N)$ on $D_1\ttimes D_k$.

Indeed, since $\phi_j:D_j\too\wdtl D_j$ is the envelope of holomorphy, for each $g_j\in\CO(D_j)$,
there exists a $\wdtl g_j\in\CO(\wdtl D_j)$ such that $g_j\equiv\wdtl g_j\circ\phi_j$. In particular,
if $\phi(z_j)=\phi(w_j)$, then $g_j(z_j)=g_j(w_j)$. Take arbitrary $c_j\in A_j$, $j=1,\dots,k$.
Then
\begin{multline*}
f(c_1,\dots,c_k,a_{k+1},\dots,a_N)=f(c_1,\dots,c_k,b_{k+1},a_{k+2},\dots,a_N)\\
=\dots=f(c_1,\dots,c_k,b_{k+1},\dots,b_N).
\end{multline*}
Thus $f(\cdot,a_{k+1},\dots,a_N)=f(\cdot,b_{k+1},\dots,b_N)$ on $A_1\ttimes A_k$. It remains to
use the identity principle.

Recall that
$$
(\phi_1\ttimes\phi_k):D_1\ttimes D_k\too\wdtl D_1\ttimes\wdtl D_k
$$
is the envelope of holomorphy (cf.~\cite{JarPfl2000}, Proposition 1.8.15 (b)). Consequently,
the function
$$
\wdtl f_\alpha(\cdot,\phi_{k+1}(a_{k+1}),\dots,\phi_N(a_N)):
=((\phi_1\ttimes\phi_k)^\ast)^{-1}(f(\cdot,a_{k+1},\dots,a_N))
$$
is well defined on
$$
\wdtl{\EX}_\alpha:=\wdtl D_1\ttimes\wdtl D_k\times\wdtl A_{k+1}\ttimes\wdtl A_N
$$
with $\wdtl f_\alpha\circ\phi=f$ on $\EX_\alpha$ and
$\sup_{\wdtl{\EX}_\alpha}|\wdtl f_\alpha|=\sup_{\EX_\alpha}|f|$.

In particular, $\wdtl f_\alpha\circ\phi=f=\wdtl f_\beta\circ\phi$ on $A_1\ttimes A_N$. Hence, by the identity
principle, $\wdtl f_\alpha=\wdtl f_\beta$ on $\wdtl{\EX}_\alpha\cap\wdtl{\EX}_\beta$.

\halfskip

Thus, we may replace $((D_j,A_j)_{j=1}^N,\BX_{N,k},\wdht{\wdtl\BX}_{N,k})$ by
$((\wdtl D_j,\wdtl A_j)_{j=1}^N, \wdtl\BX_{N,k}, \wdht{\wdtl\BX}_{N,k})$ and we may assume
$D_j$ is a domain of holomorphy and $\phi_j=\id$, $j=1,\dots,N$.

Moreover, by Remark \ref{RemOtCr}(\ref{RemOtCrAppr}), we may assume that
$A_j\subset\subset D_j\subset\subset G_j$, where $G_j$ is a Riemann domain over $\CC^{n_j}$,
$j=1,\dots,N$.

The case $k=N$ is trivial. The case $k=1$ is the classical cross theorem (Theorem \ref{ThmClassical}).
In particular,
there is nothing to prove for $N=2$. We apply induction on $N$. Suppose that the result
is true for $N-1\geq2$.

Now, we apply finite induction on $k$. The case $k=1$ is known.
Suppose that the result is true for $k-1$ with $2\leq k\leq N-1$.

Fix an $f\in\CO_s(\BX_{N,k})$ and let $C:=\sup_{\BX_{N,k}}|f|$. Recall that
$$
\BX_{N,k}=(\BX_{N-1,k-1}\times D_N)\cup(\BX_{N-1,k}\times A_N).
$$
For each $z_N\in D_N$ the function $f(\cdot,z_N)$ belongs to $\CO_s(\BX_{N-1, k-1})$.
By the inductive assumption there exists a $g_{z_N}\in\CO(\wdht{\BX}_{N-1,k-1})$
such that $g_{z_N}=f(\cdot,z_N)$ on $\BX_{N-1,k-1}$ and $\sup_{\wdht{\BX}_{N-1,k-1}}|g_{z_N}|\leq
C$. Analogously, for each $z_N\in A_N$
there exists an $h_{z_N}\in\CO(\wdht{\BX}_{N-1,k})$ such that $h_{z_N}=f(\cdot,z_N)$ on
$\BX_{N-1,k}$ and $\sup_{\wdht{\BX}_{N-1,k}}|h_{z_N}|\leq C$.
Recall that $\wdht{\BX}_{N-1,k-1}\subset\wdht{\BX}_{N-1,k}$ and
$A_1\ttimes A_{N-1}\subset\BX_{N-1,k-1}\cap\BX_{N-1,k}$. Since the set $A_1\ttimes A_{N-1}$
is not pluripolar, we get $g_{z_N}=h_{z_N}$ on $\wdht{\BX}_{N-1,k-1}$ for $z_N\in A_N$.

Consider the $2$--fold cross
$$
\BY:=\XX(\wdht{\BX}_{N-1,k-1}, A_N; \wdht{\BX}_{N-1,k},D_N)=
(\wdht{\BX}_{N-1,k-1}\times D_N)\cup(\wdht{\BX}_{N-1,k}\times A_N)
$$
and let $F:\BY\too\CC$,
$$
F(z',z_N):=\begin{cases}g_{z_N}(z'), & \text{ if }
(z',z_N)\in\wdht{\BX}_{N-1,k-1}\times D_N\\
h_{z_N}(z'), & \text{ if } (z',z_N)\in\wdht{\BX}_{N-1,k}\times A_N.
\end{cases}.
$$
Obviously, $\sup_{\BY}|F|\leq C$. To see that $F\in\CO_s(\BY)$, we have to prove that for each
$z'\in\wdht{\BX}_{N-1,k-1}$, the function $D_N\ni z_N\tuu F(z',z_N)$ is holomorphic.
We know that $F(\cdot,z_N)$ is holomorphic for each $z_N\in D_N$. Let
$$
\BZ_{N-1,k-1}:=\XX_{N-1,k-1}((A_j,D_j)_{j=2}^N).
$$
Analogously as above, for each $z_1\in D_1$
there exists a $\phi_{z_1}\in\CO(\wdht{\BZ}_{N-1,k-1})$ such that $\phi_{z_1}=f(z_1,\cdot)$
on $\BZ_{N-1,k-1}$. Thus
\begin{multline*}
F(z_1,\dots,z_N)=f(z_1,\dots,z_N)=\phi_{z_1}(z_2,\dots,z_N),\\
(z_1,\dots,z_N)\in(\BX_{N-1,k-1}\times D_N)\cap(D_1\times\BZ_{N-1,k-1})\supset
A_1\ttimes A_{N-1}\times D_N.
\end{multline*}
Consequently, $F(z',\cdot)\in\CO(D_N)$ for $z'\in A_1\ttimes A_{N-1}$ and
hence, using Terada's theorem (cf.~e.g.~\cite{Pfl2003}), we conclude that
$F\in\CO(\wdht{\BX}_{N-1,k-1}\times D_N)$.

Now, by the classical cross theorem (Theorem \ref{ThmClassical}) with $N=2$,
there exists an $\wdht f\in\CO(\wdht{\BY})$ such that
$\wdht f=F$ on $\BY$ (in particular, $\wdht f=f$ on $\BX_{N,k}$) and
$\sup_{\wdht{\BY}}|\wdht f|\leq C$. Recall that
$$
\wdht{\BY}=\{(z',z_N)\in\wdht{\BX}_{N-1,k}\times D_N:
h_{\wdht{\BX}_{N-1,k-1}, \wdht{\BX}_{N-1,k}}^\ast(z')+h_{A_N,D_N}^\ast(z_N)<1\}.
$$
Thus, to remains to apply Lemma \ref{LemFormula}.

The proof of Theorem \ref{ThmNew} is completed.


\def\bibname{References}
\bibliographystyle{amsplain}

\end{document}